\documentclass[12pt]{article}
  
\usepackage{xcolor}
\usepackage{amsmath,amssymb,amsthm}
\usepackage[framemethod=tikz]{mdframed}

\usepackage{algpseudocode,algorithm,algorithmicx}

\algrenewcommand\algorithmicrequire{\textbf{Precondition:}}
\algrenewcommand\algorithmicensure{\textbf{Postcondition:}}

\newtheorem{theorem}{Theorem}
\newtheorem{lemma}[theorem]{Lemma}
\newtheorem{problem}{Problem}
\newtheorem{conditions}{Conditions}

\title{Auxetic deformations and elliptic curves}
\author{Ciprian S. Borcea and Ileana Streinu \footnote{Rider University, NJ and Smith College, MA.}}

\begin{document}

\maketitle

\begin{abstract}
The problem of detecting auxetic behavior, originating in materials science and mathematical crystallography, refers to the property of a flexible periodic bar-and-joint framework to widen, rather than shrink, when stretched in some direction. The only known algorithmic solution for detecting infinitesimal auxeticity is based on the rather heavy machinery of fixed-dimension semi-definite programming. In this paper we present a new, simpler algorithmic approach which is applicable to a natural family of 3D periodic bar-and-joint frameworks with 3 degrees-of-freedom. This class includes most zeolite structures, which are important for applications in computational materials science. We show that the existence of auxetic deformations is related to properties of an associated elliptic curve. A fast algorithm for recognizing auxetic capabilities is obtained via the classical Aronhold invariants of the cubic form defining the curve. 
 \end{abstract}




\section{Introduction}
In this paper we address a geometric problem originating in materials science. The main result is an efficient algorithm for detecting {\em auxetic behavior} in generic 3D periodic bar-and-joint frameworks with independent edge constraints and three degrees of freedom. The method has immediate applications to frameworks of vertex-kissing tetrahedra, which include the important family of crystalline materials known as zeolites. 

\begin{figure}[h]
\centering
{\includegraphics[width=0.4\textwidth]{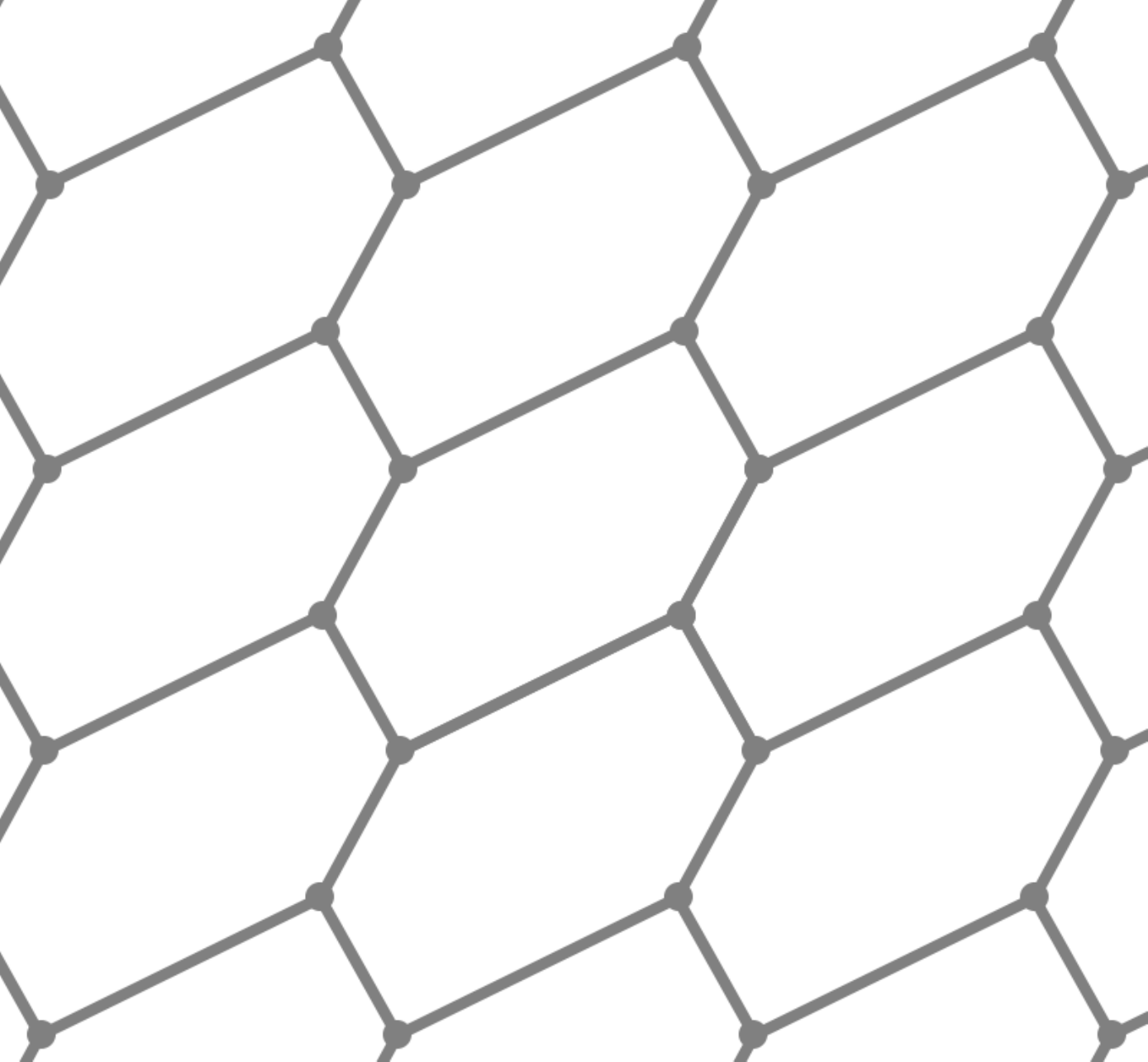}} \hspace{0.2in}
{\includegraphics[width=0.4\textwidth]{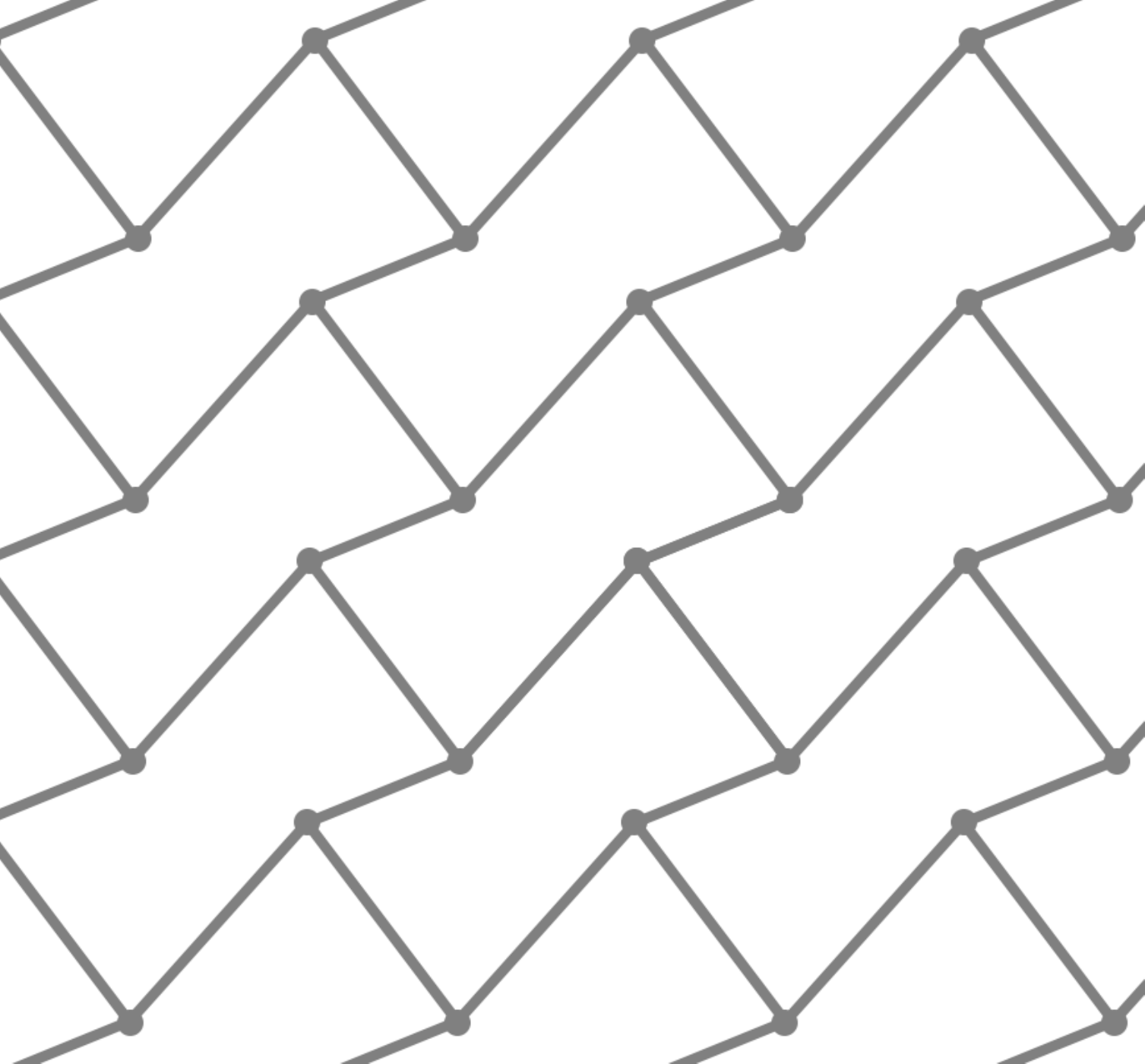}}
\caption{(Left) The {\em regular honeycomb} has positive Poisson's ratios, while the {\em re-entrant honeycomb} (right) exhibits auxetic behavior.}
\label{fig:honeycomb}
\end{figure}

\medskip\noindent{\bf Motivation: auxetic behavior in materials science.}
When a physical material is stretched along an axis, its ``typical'' response is a lateral shrinking.  Yet this behavior is not universal: certain natural or man-made materials have the rather counter-intuitive property of widening rather than shrinking. This type of ``growth'' has been called  {\em auxetic behavior} \cite{evans:etAl:molecularNetwork:Nature:1991}.  In elasticity theory, Poisson's ratio for two orthogonal directions $\mathbf{a}$ and $\mathbf{b}$ is defined as the signed ratio of the lateral effect along $\mathbf{b}$ to the extension along $\mathbf{a}$ due to uniaxial tension applied in that direction. For ``typical'' materials the Poisson's ratio is taken with a positive sign, while materials with auxetic behavior are said to have negative Poisson's ratios \cite{lakes:negativePoisson:1987,greaves:surveyPoissonRatios:resNotesRoyalSoc:2013,greaves:lakes:etAl:PoissonRatio:2011,huang:chen:negativePoisson:2016}. The difference is often illustrated in the literature with the 2D periodic bar-and-joint {\em honeycomb} examples shown in Fig.~\ref{fig:honeycomb}.  

\medskip
In mineralogy, the ``rigid unit mode'' approach to displacive phase transitions in crystalline materials relies on the periodic framework structure at the atomic level \cite{dove:displacive:1997}. In this context, auxetic behavior was observed for cristobalite \cite{yeganeh:etAl:elasticityCristobalite:1992} and conjectured among zeolites \cite{grima:etAl:doZeolitesNegativePoission:advMaterials:2000,siddornEtAl:systematicTopology:2015}. 
For cellular (periodic) materials, the importance of the underlying geometry has been frequently emphasized \cite{lee:singer:thomas:microNanoMaterials:advMat:2012,mitschkeEtAl:auxetic:rspa:2013}. However, the approach via determinations of Poisson's ratios has produced only a limited catalog of structures capable of auxetic behaviour \cite{elipe:lantada:auxeticGeometries:2012}. Currently, with new and augmented possibilities for digital fabrication (3D printing) and an increased need for synthetic {\em smart} materials, geometric structural design has reached critical importance \cite{reisEtAl:designerMatter:2015}.

\medskip
In our recent work \cite{borcea:streinu:geometricAuxetics:RSPA:arxiv:2015}, we laid the foundations of a strictly geometric theory of auxetic deformations for periodic bar-and-joint frameworks. Auxetic behavior is recognized based on the evolution of the periodicity lattice along one-parameter trajectories. As reviewed after this
introduction, the hallmark of an auxetic trajectory is that the Gram matrix of a basis of periods traces a curve of symmetric matrices with all tangent directions contained in the positive semidefinite cone. 
This purely geometric approach has uncovered new principles for auxetic design \cite{borcea:streinu:NewPrinciplesAuxeticDesign:arxiv:2016} and has opened the door to developing rigorous mathematical techniques for studying auxeticity.

\medskip\noindent{\bf Main theoretical result: overview.} 
In this paper, we investigate the auxetic capabilities of a natural class of three-dimensional periodic frameworks, which includes many structures of interest in mineralogy and materials science, such as silica and zeolites \cite{megaw:crystalStructures:1973}.  Specifically, we address the following problem. 

\medskip
\mdfsetup{roundcorner = 8pt}
\begin{mdframed}[backgroundcolor=gray!8,roundcorner=4pt]
\begin{problem}[{\bf Identify auxetic capability}]
Let ${\cal F}$ be a 3D periodic bar-and-joint frameworks with independent edge constraints, $n$ vertex orbits, $m=3n$ edge orbits and hence three degrees of freedom. Decide if ${\cal F}$ allows strictly auxetic local trajectories and if so, produce a local (infinitesimal) auxetic deformation.
\end{problem}
\ \ 
\end{mdframed}

\medskip
As it is often the case with similar questions in rigidity theory and kinematics, we assume certain genericity conditions on the framework (they will be stated explicitly in Section \ref{sec:frameworks}). In practice, most tetrahedral crystal frameworks would satisfy our conditions, which guarantee that local deformations are represented by variations of the periodicity lattice and the three dimensional tangent space at the initial configuration is in general position. 

The connection with elliptic curves indicated in the title comes from restricting the determinant function on symmetric $3\times 3$ matrices to such a three-dimensional vector subspace. The genericity condition will imply that the resulting ternary cubic form defines a non-singular projective curve. The most remarkable theoretical fact is that our question about auxetic capabilities of the periodic framework turns into a question about the real points of an elliptic curve.

\medskip\noindent{\bf Main algorithmic result: overview.} Our formulation of auxetic behavior leads to significant connections with convex algebraic geometry,  spectrahedra and semidefinite programming, as already noted in \cite{borcea:streinu:geometricAuxetics:RSPA:arxiv:2015}, Corollary 7.1. In particular, deciding infinitesimal auxeticity for a $d$-periodic bar-and-joint framework can be formulated as a semi-definite program (SDP) in fixed-dimension $d$, for which polynomial time solutions exist \cite{porkolab:khachian:semidefinite:1997}\footnote{However, as pointed out to us by Anthony Man-Cho So, this remains a purely theoretical result: no implementation of the algorithm \cite{porkolab:khachian:semidefinite:1997} is available and in practice one would have to rely on some general SDP code, such as the CVX package in MATLAB.}; see also \cite{lemon:manchoso:ye:LowRankSDP:foundTrends:2016} for a recent survey on the related topic of low-rank SDP.
Using the new connections with elliptic curves presented in this paper, we can avoid the SDP machinery for the family of $n$-vertex, $3n$-edge orbit periodic frameworks. Specifically, classical results on the Hesse form of a cubic, Aronhold invariants and the modular $J$-invariant \cite{bonifant:milnor:cubicCurves:arxiv:2016,dolgachev:classicalAlgGeom:2012,sturmfels:algorithmsInvariantTheory:2008} lead to a simpler and cleaner decision algorithm. The main steps of the algorithm involve solving an under-constrained linear system with $3n$ linear equations in $3n+3$ unknowns, along with several algebraic calculations on cubic polynomials.

\medskip\noindent{\bf Related Work.} In \cite{borcea:streinu:liftingsStresses:dcg:arxiv:2015}, we defined and studied the stronger concept of {\em expansive deformation path} for a periodic framework, and proved that it is necessarily an auxetic path. For 2D frameworks, we gave a complete combinatorial characterization (based on periodic pseudo-triangulations) of those frameworks that admit expansive deformations \cite{borcea:streinu:liftingsStresses:dcg:arxiv:2015, borcea:streinu:kinematicsExpansive:ark14:2014}. This leads, in particular, to an infinite family of planar periodic “auxetic mechanisms”. For arbitrary
dimensions, our design principles \cite{borcea:streinu:NewPrinciplesAuxeticDesign:arxiv:2016} and connections with
convex algebraic geometry bring new significance to classification problems for
spectrahedra, as pursued in \cite{ottem:etAl:quarticSpectrahedra:2015}. Also related are classical results on ternary quartics forms \cite{sottile:ternaryQuartics:2004}.

\begin{figure}[h]
\centering
{\includegraphics[width=0.4\textwidth]{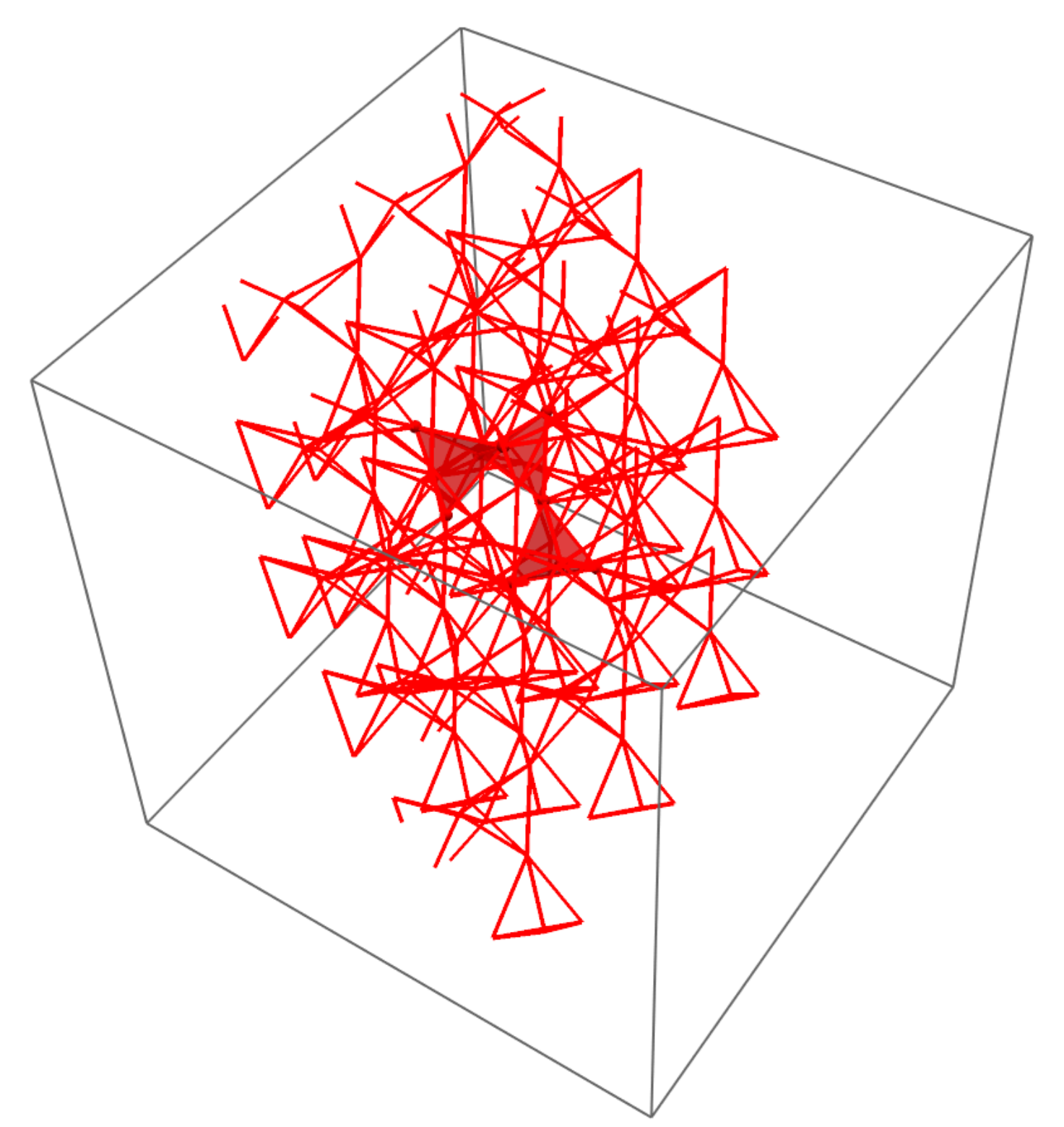}} \hspace{0.2in}
{\includegraphics[width=0.4\textwidth]{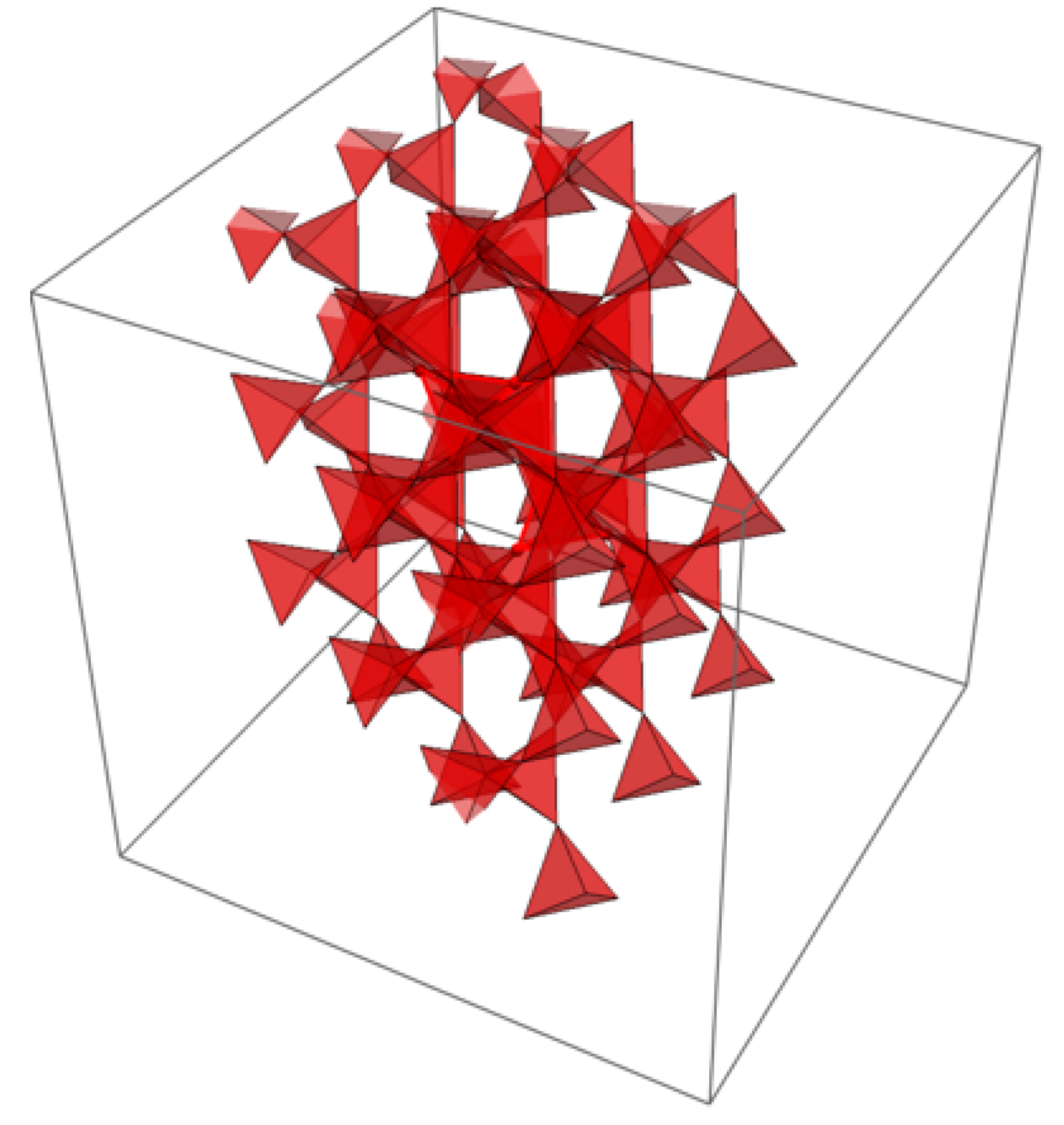}}
		\caption{The quartz framework is an example of a tetrahedral 3D periodic structure. It has $n=6$ vertex orbits and $m=3n=18$ edge orbits. The bar-and-joint framework is shown on the left. The edges of the three highlighted tetrahedra provide the 18 representatives for the edge orbits. In  \cite{borcea:streinu:geometricAuxetics:RSPA:arxiv:2015}, we show that the path of the phase transition from low to high quartz is auxetic.}
		\label{fig:quartz}
\end{figure}

\medskip\noindent{\bf Frameworks with $m=3n$ and the vertex-kissing tetrahedral structures.} The 3D frameworks studied in this paper have $n$ vertex orbits and $m=3n$ edge orbits. We assume independent edge constraints and this implies a  smooth three-dimensional local deformation space, i.e. three degrees of freedom. This choice is not accidental. In general, for arbitrary dimension $d$, the class of   frameworks with $n$ vertex orbits and $m=dn$ edge orbits has distinctive mathematical properties; for instance, the invariance of the ratio $m/n=d$ under relaxations of periodicity. A special family in dimension $2$ is the family of pointed pseudo-triangulations. Another important family, in arbitrary dimension $d$, consists of periodic frameworks made of (1-skeleta of) ``kissing-simplices'';   every vertex is common to exactly two simplices.
The 3D case of vertex-kissing tetrahedral structures includes a wealth of crystal structures studied in mineralogy, such as silicates and zeolites. Quartz, for instance, is a silicon dioxide with oxygen atoms forming a tetrahedral framework as illustrated in Fig.~\ref{fig:quartz}. Its auxetic properties, established in \cite{borcea:streinu:geometricAuxetics:RSPA:arxiv:2015}, motivated us to seek the general decision procedure presented in this paper.

\bigskip
\section{Periodic frameworks and deformations}\label{sec:definitions}

In this section we review the necessary concepts related to periodic frameworks and deformation spaces. The main references are \cite{borcea:streinu:periodicFlexibility:2010,borcea:streinu:minimallyRigidPeriodic:BLMS:2011,borcea:streinu:frameworksCrystallographicSymmetry:philTransAndArxiv:2013}. 

\bigskip\noindent{\bf Periodic graphs, vertex and edge orbits.} 
A $d$-periodic graph is a pair $(G,\Gamma)$, where $G=(V,E)$ is a simple infinite graph with vertices $V$, edges $E$ and finite degree at every vertex, and $\Gamma \subset Aut(G)$ is a free Abelian group of automorphisms which has rank $d$, acts without fixed points and has a finite number of vertex (and hence, also edge) orbits. The group  $\Gamma$ is thus isomorphic to $Z^d$ and is called the {\em periodicity group}  of the periodic graph $G$. Its elements $\gamma \in \Gamma \simeq Z^d$ are referred to as {\em periods} of $G$. 
We assume $V$ to be connected and denote by $n=|V/\Gamma|$ the number of vertex orbits and by $m=|E/\Gamma|$ the number of edge orbits.  

\bigskip\noindent{\bf Placements and frameworks.} 
A periodic placement (or simply placement) of a $d$-periodic graph $(G,\Gamma)$ in ${\mathbb{R}}^d$ is defined by two functions:
$p:V\rightarrow {\mathbb{R}}^d$ and $\pi: \Gamma \hookrightarrow {\cal T}({\mathbb{R}}^d)$, where $p$ assigns points in ${\mathbb{R}}^d$ to the vertices $V$ of $G$ and $\pi$ is a faithful representation of the periodicity group $\Gamma$, that is, an injective homomorphism of $\Gamma$ into the group ${\cal T}({\mathbb{R}}^d)$ of translations in the Euclidean space ${\mathbb{R}}^d$, with $\pi(\Gamma)$ being a lattice of rank $d$. These two functions must satisfy the natural compatibility condition $p(\gamma v)=\pi(\gamma)(p(v))$.

\bigskip
\noindent
A placement $(p,\pi)$ which does not allow the end-points of any edge to have the same image defines a {\em $d$-periodic bar-and-joint framework} 
${\cal F}=(G,\Gamma,p,\pi)$ in ${\mathbb{R}}^d$, with edges $(u,v)\in E$ corresponding to bars (segments of fixed length) $[p(u),p(v)]$ and vertices corresponding to (spherical) joints. An example in dimension $d=3$ is given in Fig.~\ref{fig:example2verts6edges}.
Two frameworks are considered equivalent when one is obtained from the other by a Euclidean isometry. 
 
\bigskip\noindent{\bf Deformation space.} 
Given a $d$-periodic framework ${\cal F}=(G,\Gamma, p,\pi)$, the collection of all periodic placements of $(G,\Gamma)$ in ${\mathbb{R}}^d$ which maintain the lengths of all edges is called the {\em realization space} of the framework. After factoring out equivalence under Euclidean isometries, we obtain the {\em configuration space} of the framework (with the quotient topology). The {\em deformation space} is the connected component of the configuration space which contains the initial framework.

\bigskip\noindent{\bf Factoring out the trivial isometries.} 
A convenient way to factor out equivalence under Euclidean isometries is to use coordinates relative to a basis given by generators of the periodicity lattice and to retain the metric information via the Gram matrix 
of the chosen basis. The procedure is described in \cite{borcea:streinu:frameworksCrystallographicSymmetry:philTransAndArxiv:2013}, Section 4,
and will be recalled here concisely.

\bigskip
\noindent
After choosing an independent set of $d$ generators for the periodicity lattice $\Gamma$, the image $\pi(\Gamma)$ is completely described via the $d\times d$ matrix $\Lambda$ with column vectors $(\lambda_i)_{i=1,\cdots,d}$ given by the images of the generators under $\pi$. The Gram matrix for this basis will be $\omega=\Lambda^t\cdot \Lambda$.

\bigskip \noindent
Let us fix now a complete set of vertex representatives $v_0,v_1,...,v_{n-1}$
for the $n$ vertex orbits of $(G,\Gamma)$. The framework ${\cal F}$ has them
positioned at $p_i=p(v_i)$. When we pass from these Cartesian coordinates to
lattice coordinates $q_i$, we consider $v_0$ to be the origin, that is $q_0=0$,
and then $\Lambda q_i=p_i-p_0$.

\bigskip \noindent
In this manner, the Cartesian description of the framework ${\cal F}$ given by $(p_0,...,p_{n-1})\in ({\mathbb{R}}^d)^n$ and $\Lambda \in GL(d)$, which requires $dn+d^2$
parameters, is reduced to $(q_1,...,q_{n-1},\omega)$, which involves only
$d(n-1)+{d+1\choose 2}$ parameters, since $\omega$ is a symmetric $d\times d$
matrix. The dimensional drop of ${d+1\choose 2}$ reflects the factoring out of the
action of the Euclidean group of isometries $E(d)$, which has this dimension.

\bigskip \noindent
We recall here the form of the equations expressing the constant (squared) length of edges, when using parameters $(q_1,...,q_{n-1},\omega)$ (cf. \cite{borcea:streinu:frameworksCrystallographicSymmetry:philTransAndArxiv:2013}, formula (4.1)). Let us consider an edge (denoted here simply by $e_{ij}$) which goes from $v_i$ to a vertex in the orbit of $v_j$. Then, in {\em Cartesian coordinates}, the edge vector is given by 
$p_j+\lambda_{ij}-p_i$, with some period $\lambda_{ij}=\Lambda n_{ij}\in \pi(\Gamma)$ and $n_{ij}\in Z^d$. In {\em lattice coordinates}, the edge vector is given by $q_j+n_{ij}-q_i$ and the squared-length equation is:
\begin{align}
	\ell(e_{ij})^2 =  \ & \  || p_j+\lambda_{ij}-p_i ||^2= \nonumber \\
    =\langle \omega (q_j+n_{ij}-q_i),(q_j+n_{ij}-q_i) \rangle = \ & \  
	(q_j+n_{ij}-q_i)^t \omega(q_j+n_{ij}-q_i) \label{sqLength}
\end{align}

\noindent
This is the source of the linear system for computing infinitesimal deformations, as used in subsequent sections. Besides the elimination of equivalence under isometries, an obvious advantage of this formulation is that the Gram matrix  $\omega=\Lambda^t\cdot \Lambda$ of the chosen basis of periodicity generators appears explicitly. 

\bigskip

\section{Auxetics and spectrahedra}\label{spectrahedra}

We review now, following \cite{borcea:streinu:geometricAuxetics:RSPA:arxiv:2015}, the fundamental notions of geometric auxetics. 

\bigskip\noindent{\bf One-parameter periodic deformations.} 
A {\em one-parameter deformation of the periodic framework} ${\cal F}=(G,\Gamma, p,\pi)$ is a smooth family of placements  $p_{\tau}: V\rightarrow {\mathbb{R}}^d$,  parametrized by time $\tau \in (-\epsilon, \epsilon)$ in a small neighborhood of the initial placement $p_0=p$, which satisfies two conditions: (a) it maintains the lengths of all the edges $e\in E$, and (b) it maintains periodicity under $\Gamma$, via faithful representations $\pi_{\tau}:\Gamma \rightarrow {\textbf T}({\mathbb{R}}^d)$ which {\em may change with $\tau$ and give an associated variation of the periodicity lattice of translations} $\pi_{\tau}(\Gamma)$. 

\bigskip \noindent
With our chosen generators for the periodicity group $\Gamma$, we have, at any
moment of time $\tau \in  (-\epsilon, \epsilon)$, a lattice basis $\Lambda_{\tau}$ and the corresponding Gram matrix
$\omega_{\tau}=\omega(\tau)=\Lambda^t_{\tau}\Lambda_{\tau}$.

\bigskip\noindent{\bf Auxetic path.} A one-parameter deformation  $(G,\Gamma, p_{\tau},\pi_{\tau}), \tau\in (-\epsilon,\epsilon)$ is called an {\em auxetic path}, or simply {\em auxetic}, when the curve of Gram matrices $\omega(\tau)$ has all its
tangents in the cone of positive semidefinite symmetric $d\times d$ matrices. When all tangents are in the positive definite cone, the deformation is called {\em strictly auxetic.}

\medskip \noindent
When we operate in lattice coordinates $(q_1,...,q_{n-1},\omega)$, we use the abbreviations $\dot{q_i}=\frac{d q_i}{d \tau}(0)$ and  $\dot{\omega}=\frac{d \omega}{d \tau}(0)$. From (\ref{sqLength}) we obtain
the linear system for infinitesimal deformations $(\dot{q}_1,...,\dot{q}_{n-1},\dot{\omega})$:
\begin{align}
& \langle \dot{\omega} e_{ij}, e_{ij} \rangle + 2\langle \omega e_{ij}, \dot{e}_{ij}\rangle = 0 
\label{eq:infFlexes}
\end{align}

\noindent
where $e_{ij}=q_j+n_{ij}-q_i$, with $n_{ij}\in Z^d$, is meant to run over a complete set of representatives of edge orbits and $\dot{e}_{ij}=\dot{q}_j-\dot{q}_i$.  For notational simplicity, the edge symbol $e_{ij}$ omits the specific period shift, but it will be understood that several  such edge representatives (with various shifts $n_{ij}$) may be implicated in the system. There are $d(n-1)+{d+1\choose 2}=dn +{d\choose 2}$ unknowns and $m=|E/\Gamma|$ equations.

\bigskip \noindent
Infinitesimally, our auxetic criterion requires $\dot{\omega}$ to be a positive
semi-definite matrix, while strict infinitesimal auxeticity needs a positive
definite $\dot{\omega}$. It follows immediately that {\em infinitesimal auxetic deformations} define a {\em cone} in the tangent space at ${\cal F}$ of the framework's deformation space. We speak of {\em infinitesimal auxetic capabilities} when there are infinitesimal deformations with $\dot{\omega}\neq 0$ and $\dot{\omega}$ positive semidefinite.

\bigskip \noindent
In convex optimization theory and semidefinite programming, an affine section
of the positive semidefinite cone in the linear space of symmetric $d\times d$
matrices is called a {\em spectrahedron} \cite{vinzant:whatIsSpectrahedron:2014,ottem:etAl:quarticSpectrahedra:2015}. When the section goes through
the origin we have a {\em spectrahedral cone}. Thus, when we look for infinitesimal
auxetic capabilities, we have to look at the image of the space of infinitesimal
deformations in the space of symmetric $d\times d$ matrices and examine its intersection with the positive semidefinite cone. 

\bigskip \noindent
This important connection between infinitesimal auxetic capabilities of a periodic
framework and semidefinite programming was established in \cite{borcea:streinu:geometricAuxetics:RSPA:arxiv:2015}, Corollary 7.1
and will be a pivotal element of the current investigation, which is targeted on
a specific class of three-dimensional periodic frameworks with three degrees of freedom. From this perspective, our inquiry will be concerned with the six-dimensional vector space of $3\times 3$ symmetric matrices and spectrahedral
cones resulting from the intersection of the positive semidefinite cone with
(generic) three-dimensional vector subspaces.

\bigskip
\section{Regular three-dimensional periodic frameworks}\label{sec:frameworks}

There are several motivations for looking at three-dimensional periodic frameworks with three degrees of freedom. Applications to crystalline materials evidently require dimension three and important families studied for properties related to framework flexibility have, geometrically, three degrees of freedom. Zeolites, for instance, have framework structures made of vertex sharing tetrahedra and generic periodic frameworks of this type have, locally, a three-dimensional deformation space \cite{borcea:streinu:periodicFlexibility:2010}. 

\bigskip \noindent
If we assume independent edge constraints in the linear system \eqref{eq:infFlexes} obtained above for infinitesimal deformations, we see that a three dimensional solution space amounts to $m=3n$. Mathematically, this type of periodic framework is quite relevant for periodicity relaxation problems, since this proportion of edge orbits to vertex orbits remains the same when adopting a sublattice of periods (of finite index).

\bigskip \noindent
As with similar definitions in rigidity theory, we have to assume certain genericity conditions on the framework, in particular independence of the edges \cite{borcea:streinu:minimallyRigidPeriodic:BLMS:2011}. For reasons that will become clear along the way, in this paper we consider only {\em regular periodic frameworks} in ${\mathbb{R}}^3$, defined by the following: 

\bigskip

\mdfsetup{roundcorner = 8pt}
\begin{mdframed}[backgroundcolor=gray!8,roundcorner=4pt]
\begin{conditions}[{\bf Regularity conditions}]
\ \ \ 
\begin{enumerate}
\item[i.] The number of edge representatives is $m=3n$.
\item[ii.] The edge constraints are independent, that is, the system \eqref{eq:infFlexes} has maximal rank $m$.
\item[iii.] The projection from infinitesimal deformations to $3\times 3$ symmetric matrices:
$$(\dot{q}_1,...,\dot{q}_{n-1},\dot{\omega})\mapsto \dot{\omega}$$
is one-to-one.
\item[iv.] The projective cubic curve obtained by restricting the determinant function to the image of the above projection is non-singular.
\end{enumerate}
\end{conditions}
\ \ 
\end{mdframed}

We can now formulate precisely the problem solved in this paper:

\medskip
\mdfsetup{roundcorner = 8pt}
\begin{mdframed}[backgroundcolor=gray!8,roundcorner=4pt]
\begin{problem}[{\bf Identify auxetic capability}]
Let ${\cal F}$ be a 3D {\em regular} periodic bar-and-joint framework. Decide if ${\cal F}$ allows strictly auxetic deformations and if so, produce a local (infinitesimal) auxetic deformation.
\end{problem}
\ \ 
\end{mdframed}
\medskip\noindent{\bf Remarks.}\ It follows from the implicit function theorem, under the regularity conditions stated above,  that the local deformation space of ${\cal F}$ is a {\em smooth threefold}. This means, in particular, that it is enough to solve the infinitesimal problem of {\em deciding if ${\cal F}$ has strictly auxetic infinitesimal deformations}, since any curve of symmetric matrices with a positive definite tangent directions at one point will maintain this property in a neighborhood of that point. In other words, ${\cal F}$ {\em allows strictly auxetic deformations} if and only if it {\em allows strictly auxetic infinitesimal deformations.}

\medskip \noindent
Moreover, under the genericity condition (iv), intersecting the positive definite cone is equivalent with intersecting non-trivially its closure, the positive semidefinite cone,
since we have a transversal intersection with the boundary in the generic case.
Thus, our decision problem turns into the general question about infinitesimal
auxetic capabilities which can be solved by semidefinite programming as stated
in \cite{borcea:streinu:geometricAuxetics:RSPA:arxiv:2015}, Corollary 7.1. The solution presented in this paper for our distinctive class of frameworks shows that geometric characteristics can be used for alternative algorithms, and thus avoid the general machinery of semidefinite programming \cite{porkolab:khachian:semidefinite:1997}.

\section{Ternary cubic forms and their invariants}\label{sec:cubics}

In this section we explain our geometrical approach, present the theoretical basis of the decision algorithm and review the classical results about elliptic curves. 

\medskip \noindent
We have shown above that the question about strictly auxetic deformations of
a regular periodic framework ${\cal F}$ leads to an intersection problem in
the space of $3\times 3$ symmetric matrices between the positive semidefinite 
cone and a three-dimensional vector subspace image of the infinitesimal deformations of ${\cal F}$. Strictly auxetic deformations exist if and only if
the three-dimensional vector subspace cuts through the positive definite cone.

\medskip \noindent
Our approach is to investigate the resulting intersection with the boundary of the positive definite cone, which is contained in the vanishing locus of the determinant function. Since this function is a homogeneous cubic form, its restriction to our three-dimensional subspace gives a {\em ternary cubic form}
and projectively, a real cubic curve in $P_2({\mathbb{R}})$. The genericity condition (iv)
allows the assumption that this cubic curve is non-singular (over the complex numbers ${\mathbb{C}}$).

\medskip \noindent
Non-singular projective cubics are also known as {\em elliptic curves} and there is a wealth of classical results about them. We limit the main references to 
\cite{bonifant:milnor:cubicCurves:arxiv:2016,dolgachev:classicalAlgGeom:2012,sturmfels:algorithmsInvariantTheory:2008}, which cover all the facts needed for our arguments.  Different
sources may have different choices of coefficients in front of certain invariants.
Since we are concerned with real ternary cubics, we follow mostly
\cite{bonifant:milnor:cubicCurves:arxiv:2016}.

\medskip \noindent
We consider a three-dimensional real vector space with coordinates $(x,y,z)$.
We are interested in real cubic forms $f(x,y,z)$ and the projective cubic curves
$f(x,y,z)=0$ they define in $P_2({\mathbb{R}})$, with projective coordinates $(x:y:z)$.
There are ten distinct monomials of degree three in $x,y,z$ and a ternary cubic form is determined by the ten coefficients of these monomials. 

\medskip\noindent{\bf Hesse normal form.} A non-singular real projective cubic $f(x,y,z)=0$
is projectively equivalent over ${\mathbb{R}}$ with a unique cubic of the form:
\begin{align}\label{eq:Hesse}
&	 x^3 + y^3 + z^3 - 3k xyz =0 
\end{align}	
	
\noindent
with $k\in {\mathbb{R}}, \ k\neq 1$. This form is called the Hesse normal form of the cubic. 
The three real inflection points are on the line $x+y+z=0$, at $(0:1:-1)$, $(-1:0:1)$ and $(1:-1:0)$.

\medskip\noindent{\bf Aronhold invariants $S$ and $T$.} The action of the special linear group
$SL(3)$ on the coordinates $(x,y,z)$ induces a linear representation on the ten-dimensional vector space of ternary cubics. Homogeneous polynomials in the 
coefficients of ternary cubics are called invariant when they remain unchanged
under this induced action. The Aronhold invariants $S$ and $T$ are given by two explicit invariant polynomials (with rational coefficients) of degree four, respectively six \cite{sturmfels:algorithmsInvariantTheory:2008}. 

\medskip\noindent{\bf Discriminant.}\ The discriminant $\Delta$ is the invariant polynomial of degree twelve given by:
\begin{align}\label{Delta}
&\Delta = (4S)^3+T^2
\end{align}

\noindent
A cubic form defines a non-singular curve in $P_2({\mathbb{C}})$ if and only if 
$\Delta\neq 0$. For negative $\Delta$, the real cubic in $P_2({\mathbb{R}})$ is connected,
while for positive $\Delta$, it has two connected components.

\medskip\noindent{\bf Modular invariant $J$.}\ For non-singular cubic forms, the {\em modular invariant} $J$ is defined as a rational expression in the Aronhold invariants $S$
and $T$ by the formula:
\begin{align}\label{eq:modulus}
&J=\frac{(4S)^3}{(4S)^3+T^2}
\end{align}

\noindent
Since both numerator and denominator are homogeneous polynomials of degree twelve,
the modular invariant $J$ is invariant under the general linear group $GL(3)$.
Moreover, two non-singular cubic curves in $P_2({\mathbb{C}})$ are projectively equivalent
over ${\mathbb{C}}$ if and only if $J$ gives the same value when evaluated for the two cubic forms in the defining equations. This value is called the {\em modulus} of the curve.

\medskip \noindent
For a cubic in {\em Hesse normal form} \eqref{eq:Hesse}, we refer to the corresponding ternary form directly through the Hesse parameter $k$. Then, the values $S(k)$ and $T(k)$ are related to the $k$ by the following equations \cite{dolgachev:classicalAlgGeom:2012}.
\begin{align}\label{eq:S}
& S=-\frac{1}{2}k-\frac{1}{2^4}k^4
\end{align}
\begin{align}\label{eq:T}
& T=1+\frac{5}{2}k^3-\frac{1}{2^3}k^6
\end{align}

\noindent
The discriminant takes the form:
\begin{align}\label{eq:Delta_k}
& \Delta(k) = (4S(k))^3+T(k)^2 = (1-k^3)^3
\end{align}

\noindent
and the case $k=1$ in \eqref{eq:Hesse} is detected as singular.

\medskip \noindent
For the modular invariant we find:
\begin{align}\label{eq:J}
& J(k) =\frac{(4S(k))^3}{(4S(k))^3+T(k)^2} = \frac{k^3(k^3+8)^3}{64(k^3-1)^3}
\end{align}

\noindent
and the right hand side is invariant under the involution $\eta(k)=\frac{k+2}{k-1}$. This means that $J(\eta(k))=J(k)$, and the two curves are projectively equivalent over the complex numbers ${\mathbb{C}}$. The two fixed points of the involution $\eta$
are $k_{\pm}=1\pm \sqrt{3}$ and $J(k_{+})=J(k_{-})=1$. Thus, for {\em real Hesse cubics} i.e. $k\in {\mathbb{R}}, k\neq 1$, for any given real modulus there are exactly two real Hesse parameters with that modulus. The two curves {\em are not} projectively equivalent over ${\mathbb{R}}$, but as complex curves, they are projectively equivalent over ${\mathbb{C}}$. For more details and illustrations, we refer to \cite{bonifant:milnor:cubicCurves:arxiv:2016}.

\medskip \noindent
For any non-singular cubic form $f$, we first compute the (values of) the Aronhold invariants $S(f)$ and $T(f)$ and then use formula \eqref{eq:modulus} to
compute the modulus $J(f)$. In the {\em real case}, we know that there is a {\em unique} real parameter $k$ which gives a Hesse curve projectively equivalent over ${\mathbb{R}}$ with the curve $f=0$, although, as reviewed above, there's another real parameter yielding the same modulus $J(f)$. The proper Hesse parameter $k$ can be
recognized as follows.

\begin{lemma}\label{find_k}
Let $f$ be a non-singular ternary cubic form with Aronhold invariants $(S(f),T(f))$ and modulus $J(f)$. Then, the quartic
equation in $k^3$ resulting from relation \eqref{eq:J}, namely:
\begin{align*}
	& 64J(f)(k^3-1)^3=k^3(k^3+8)^3
\end{align*}

\noindent
has exactly two distinct real solutions, hence two distinct real solutions $k_1$ and $k_2$. If $T(f)\neq 0$, i.e. $J(f)\neq 1$, then $T(k_1)$ and $T(k_2)$ have opposite signs and the proper solution which gives a Hesse cubic projectively equivalent over ${\mathbb{R}}$ with the curve $f=0$ must share the sign of $T(f)$. If $T(f)=0$, then $S(k_1)$ and $S(k_2)$ have opposite signs and the proper solution must give the sign of $S(f)$.
\end{lemma}

\medskip \noindent
After this brief review of relevant classical notions and results, we proceed now with our inquiry. Let us observe first that, if our (generic) three-dimensional subspace of $3\times 3$ symmetric matrices cuts through the positive definite cone, the corresponding (non-singular) real cubic curve must have {\em two connected components}. This follows from the fact that the intersection
with the boundary of the positive definite cone gives one component which appears
as a convex curve in adequately chosen affine planes. Since this component cannot
contain inflection points, there must be a second component. We obtain:

\begin{lemma}\label{lem:negative}
If the discriminant of the ternary cubic $f$ associated to the periodic framework
${\cal F}$ 	is negative, i.e. $\Delta(f) < 0$, then the framework does not allow
non-trivial auxetic deformations.
\end{lemma}
	
\medskip
However, a positive discriminant is not a sufficient condition for strictly
auxetic capabilities. We must check if the connected component without inflection points actually corresponds with positive semidefinite matrices. This can be determined by finding the explicit {\em real linear transformation} which takes the original cubic form to its Hesse normal form. Note that, by Proposition~\ref{find_k}, the Hesse parameter can be determined from the invariants $(S(f),T(f),J(f))$. Then, the decision is {\em yes} if the preimage of $(1:1:1)$
corresponds to a positive or negative definite matrix and {\em no} otherwise.

\medskip
We explain this procedure in more detail. For a cubic in Hesse normal form
\eqref{eq:Hesse}, the tangent lines at the three inflection points $(0:1:-1)$, $(-1:0:1)$ and $(1:-1:0)$ have equations $kx+y+z=0$, $x+ky+z=0$, respectively
$x+y+kz=0$. For non-singular real cubics (i.e. $k\neq 1$) these lines are concurrent only for $k=-2$, when the real curve is connected. Thus, for
real cubics with two real components (i.e. $k > 1$), the three tangents and
the line $x+y+z=0$ of the inflection points give four points in general
position in the dual projective plane $P_2({\mathbb{R}})^*$.

\medskip \noindent
It follows that, if we know the three real inflection points of a non-singular
real cubic with two connected components, we can compute the three tangent lines
and, after choosing a matching order for the inflection points, there is a unique
real projective transformation which takes the line of real inflection points to the line
$x+y+z=0$ of the cubic's Hesse normal form and matches the inflection points and
tangents. Obviously, this projective transformation takes the cubic to its
Hesse normal form. A different choice of matching order for the inflection points
amount to composition with a permutation of the variables $x,y,z$.

\medskip \noindent
An algorithm for computing the three real inflection points of a non-singular real cubic is described in \cite{chen:wang:inflectionPoints:2003}.

\medskip \noindent
For the Hesse form, the point $(1:1:1)\in P_2({\mathbb{R}})$, as the only invariant point under permutations,  must be in the
interior of the triangle formed by the three tangents in the affine plane which
has the line of inflection points at infinity. Since each connected component goes to itself under permutations, the component without inflection points must
appear as a convex curve around the invariant point (and be contained in the triangle of tangents as well). Again, illustrations may be found in \cite{bonifant:milnor:cubicCurves:arxiv:2016}.

\section{Decision algorithm for auxetic capabilities}\label{sec:decision}

A high level description of the algorithm is given below in Algorithm  \ref{alg:auxeticDecision}. 

\medskip

\algblock[Name]{Step}{EndStep}

\begin{algorithm}
  \caption{Decide existence of strictly auxetic deformations.
    \label{alg:auxeticDecision}}
  \begin{algorithmic}
	 
	 \Require{${\cal F}$ is a regular periodic framework.}
	    \Statex \Function{Auxetic}{${\cal F}$} 
				\Step \ 1.
				\State - Set up the linear system for periodic infinitesimal deformations.
				\State - Solve it in terms of $3$ independent variables chosen from the $6$ 
				\State giving the infinitesimal deformations of the Gram matrix.
				\State - If this is not possible, STOP: the framework is not regular.
				\EndStep \ 1.
				\State
				\Step \ 2.
				\State - $M \leftarrow$ Substitute the resulting linear forms in the $3\times3$ matrix of 
				\State infinitesimal deformations of the Gram matrix.
				\State - Compute the determinant $Det(M)$. The result is a cubic form 
				\State $C(X,Y,Z)$ in $3$ variables, called $X, Y$ and $Z$.
				\EndStep \ 2.
				\State
				\Step \ 3.
				\State - Compute the Aronhold invariants S and T for C(X,Y,Z).
				\State - Using S and T, compute the discriminant D of C(X,Y,Z).
				\State - If $D=0$, cubic is singular. \Return  ``not regular''.
				\State - If $D < 0$, \Return NO: the framework does not have infinitesimal 
				\State auxetic deformations.
				\EndStep \ 3.
				\State
				\Step \ 4.
				\State If $D > 0$: 
				\State - compute the $3\times3$ linear matrix $L$ for the transformation 
				\State of the cubic $C(X,Y,Z)$ to the Hesse normal form $H(x,y,z)$.  
				\State - Compute the pre-image of the point $(1,1,1)$. This is a vector of 
				\State specific values for $(X,Y,Z)$. 
				\State - The corresponding constant symmetric matrix $M$ is then tested  
				\State for being definite (either positive or negative definite) and the 
				\State corresponding output is produced: 
				\State - \Return YES, if definite or  \Return NO, if indefinite.
				\EndStep \ 4.
	\State
	\EndFunction
\end{algorithmic}
\end{algorithm}

\medskip\noindent{\bf Input and Output.} The input is a framework ${\cal F}$, with quotient graph $G=(V,E)$, $|V|=n$ and $|E|=3n$, periodicity lattice given by the Gram matrix $\omega$, and vertex representatives $q_i$ presented in ``lattice'' or ``crystallographic coordinates'', as described in Section \ref{sec:definitions}. In particular, $q_0$ is fixed at the origin, and each of the remaining $n-1$ points has 3 coordinates in the interval $[0,1)$. In a pre-processing step we can check that the edges are independent, for example by writing the $3n \times 3n$ rigidity matrix of the framework and computing its rank. The algorithm proceeds only if the rank is maximum. If the framework is not regular (cf. definition in Section \ref{sec:frameworks}), the algorithm will stop and report that it cannot give a precise answer. Otherwise, it reports the presence or not of auxetic capabilities.  We remark that the algorithm runs on all {\em generic inputs}. In particular, the set of inputs on which it will stop without giving a definitive answer has measure zero. 

\medskip\noindent{\bf Analysis.} In Step 1, we set up the $3n$ linear equations \eqref{eq:infFlexes} with $3n+3$ unknowns $(\dot{q}_1,...,\dot{q}_{n-1},\dot{\omega})$ corresponding to the infinitesimal flexes $\dot{q}_i$ of the lattice coordinates $q_i, i=1,\cdots,n-1$ of the framework and infinitesimal changes $\dot{\omega}$ in the $6$ parameters describing the symmetric Gram matrix $\omega$. We solve the linear system, and attempt to express $3n$ of the dependent variables in terms of $3$ free variables chosen from among the $\dot{\omega}$. If this is not possible, then we stop and report that the framework is not regular. In Step 2, we compute the determinant of the $3\times 3$ infinitesimal Gram matrix $\dot{\omega}$, whose entries are linear forms in the three free variables and obtain the ternary cubic described in Section \ref{sec:cubics}. In Step 3, we compute the S and T Aronhold invariants of the cubic and the discriminant $\Delta$ described in Section \ref{sec:cubics}. We stop if $\Delta=0$: the framework violates the regularity conditions. If $\Delta < 0$, the framework does not have auxetic capabilities, cf. Lemma \ref{lem:negative}. Otherwise, if $\Delta > 0$, we apply the procedure described in Section \ref{sec:cubics} after Lemma \ref{lem:negative} to distinguish, correctly, the auxetic and non-auxetic cases.

\medskip\noindent{\bf Complexity.} The running time of the algorithm is dominated by Step 1 (Gaussian elimination) and takes $O(n^3)$ arithmetic operations. The other steps are constant time calculations involving equations of degree at most 3.

\medskip\noindent{\bf Extensions.} The decision algorithm can be easily modified to return an actual auxetic infinitesimal deformation. Indeed, the pre-image of the point $(1:1:1)$ in step 4 gives specific values for the free variables $X,Y,Z$ of the cubic from Step 2, and all the infinitesimal deformation variables (for the points $q$ and for the Gram matrix $\omega$) can be expressed in terms of these three values, cf. Step 1.  Finally, the step-by-step calculation of infinitesimal auxetic deformations can be used in a standard numerical simulation, with a sufficiently small time step, of an auxetic trajectory, if one exists.

\section{Examples}\label{sec:examples}

\begin{figure}[h]
\centering
{\includegraphics[width=0.34\textwidth]{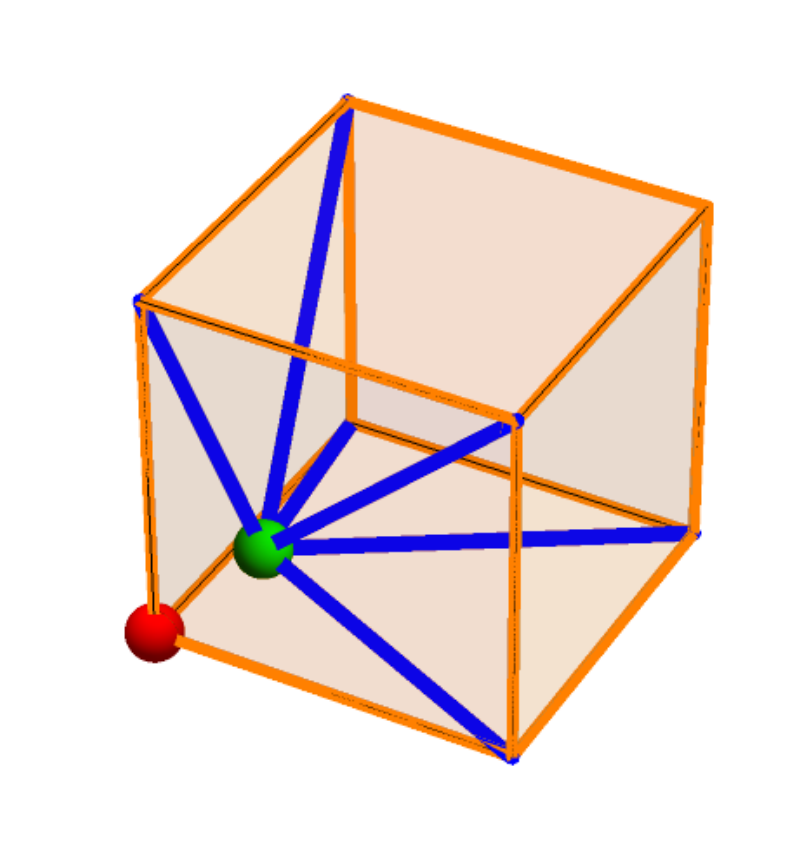}}
{\includegraphics[width=0.45\textwidth]{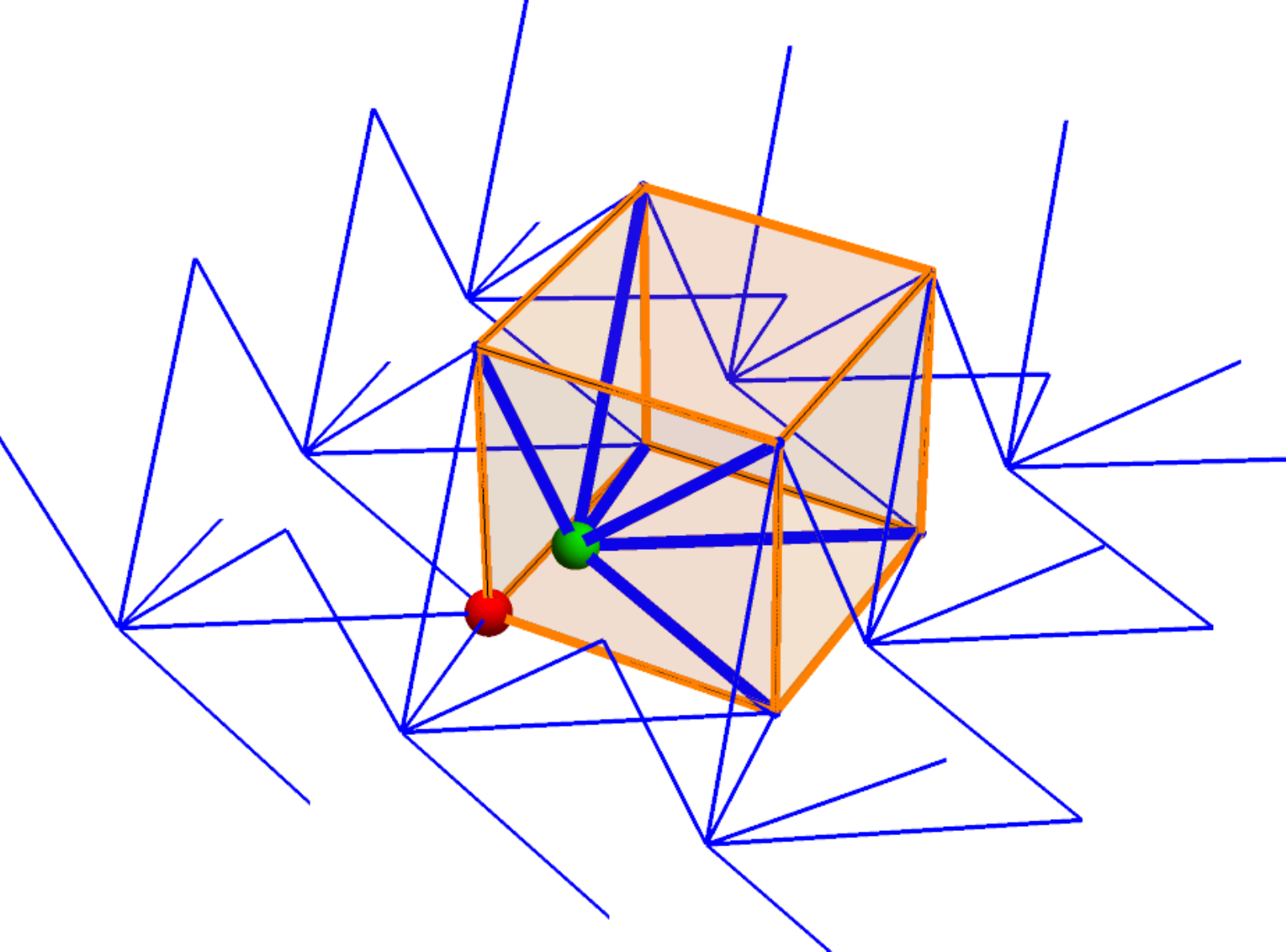}}
\caption{A 3D periodic framework with $n=2$ vertex and $m=3n=6$ edge orbits. (Left) The unit cell with the $2$ vertex representatives (colored) and the $6$ edge representatives.  (Right) A small $2\times 2\times 1$ fragment of the infinite framework.}
\label{fig:example2verts6edges}
\end{figure}

We now illustrate our method with a type of structure as shown in Fig.~\ref{fig:example2verts6edges}. We define a family of periodic frameworks ${\cal F}(\lambda)$ in ${\mathbb{R}}^3$ with $n=2$ vertex orbits and $m=3n=6$ edge orbits. The periodicity lattice is the standard integer lattice $Z^3$ and the chosen generators are the vectors $e_i, i=1,2,3$ of the standard basis. The first orbit of vertices is represented by the origin (the red point in Fig.~\ref{fig:example2verts6edges}), and coincides with $Z^3$. For ${\cal F}(\lambda)$, the second orbit of vertices is represented by the point $p=p(\lambda)=\lambda(e_1+e_2+e_3)$ (the green point in Fig.~\ref{fig:example2verts6edges}). The edge orbits are represented by six edges connecting $v$ with the vertices of the other orbit placed at $e_1,e_2,e_3, e_1+e_2, e_2+e_3, e_3+e_1$ (the thick blue edges in Fig.~\ref{fig:example2verts6edges}).

In lattice coordinates $(q,\omega)$, the framework ${\cal F}(\lambda)$ 
will have $q=\lambda(e_1+e_2+e_3)$ and $\omega=I_3$, the identity matrix.
Although $q=p$ in our initial setting, we maintain the notational distinction
between lattice coordinates and Cartesian coordinates. Thus, the six edge vectors
in lattice coordinates are written as $E_i = e_i-q, \  i=1,2,3$ and $E_{jk} = e_j+e_k-q, \ ij\in \{ 12, 23, 13\}$,
and the linear system \eqref{eq:infFlexes} for infinitesimal deformations takes the form:
\begin{align*}
\langle \dot{\omega} E_i, E_i \rangle - 2\langle E_{i}, \dot{q}\rangle = 0  \ & \ \mbox{for} \ \ i = 1,2,3 \\
\langle \dot{\omega} E_{jk}, E_{jk} \rangle - 2\langle E_{jk}, \dot{q}\rangle = 0 \  &  \ \mbox{for} \ \ j,k = 1,2,3, j\neq k
\end{align*}
	
\medskip \noindent
Exploiting permutation symmetry and choosing the diagonal elements of $\dot{\omega}$ as parameters, the three-dimensional space of solutions $\dot{\omega}$ is given by symmetric $3\times 3$ matrices with equal off-diagonal entries of the form:
\begin{align}\label{eq:3space} 
&	\dot{\omega}_{kl}= \frac{\lambda (1-\lambda)}{2-6\lambda(1-\lambda)}\cdot 
(\sum_{i=1}^3 \dot{\omega}_{ii}) 
\end{align}

\noindent
We notice the invariance of the
fractional coefficient in \eqref{eq:3space} under the involution 
$\lambda \mapsto 1-\lambda$, which comes from the isomorphism
of the frameworks ${\cal F}(\lambda)$ and ${\cal F}(1-\lambda)$ when reflecting
in the plane $x_1+x_2+x_3=1/2$ (in Cartesian coordinates).

\medskip
For simpler expressions, we introduce the following notations:

\begin{align}\label{eq:notation} 
&& \ell = \lambda (1-\lambda), \ \ \
\mu=\frac{\lambda (1-\lambda)}{2-6\lambda(1-\lambda)}=\frac{\ell}{2-6\ell},\ \ \
\rho=\mu^2(2\mu-1)
\end{align}

\medskip \noindent
We notice that, for $\lambda \in (-\infty, \infty)$, we have $\ell \in (-\infty,1/4]$ and $\mu \in (-1/6,1/2]$, with the resulting range 
$\rho \in [-1/27,0]$.

\medskip \noindent
We also rename our parameters as $X= \dot{\omega}_{11}$, $Y= \dot{\omega}_{22}$
and $Z=	\dot{\omega}_{33}$. Then the cubic equation given by $det(\dot{\omega})=0$ becomes:
\begin{align}\label{eq:Det}
& XYZ+\rho(X+Y+Z)^3=0
\end{align}

\noindent
It is easily verified that the points $(0:1:-1)$, $(-1:0:1)$ and $(1:-1:0)$
are inflection points of this cubic curve. The corresponding tangent lines 
are $X=0$, $Y=0$ and $Z=0$.

\medskip 
\noindent
The last part of the algorithm requires the computation of the real projective transformation which takes a given cubic form to the Hesse normal form:
\begin{align}\label{eq:HesseBis}
&	 x^3+y^3+z^3-3k xyz=0. 
\end{align}	

\noindent
In this case, it is easy to see that (up to permutation), the relation
between \eqref{eq:Det} and \eqref{eq:HesseBis} takes the form:
$$ X=kx+y+z $$
$$ Y=x+ky+z $$
$$ Z=x+y+kz $$

\noindent
with the parameters $\rho$ and $k$ related by:
\begin{align}\label{eq:for_k}
&	\rho=-\frac{k^2+k+1}{3(k+2)^3}
\end{align}

\noindent
This formula follows from the necessary vanishing of monomial coefficients 
for $x^2y$ etc. The singular case $k=1$ for \eqref{eq:HesseBis} corresponds with
the singular case $\rho=-1/27$ for \eqref{eq:Det}. The case $k=-2$ cannot arise,
since $\rho=\pm\infty$ is excluded.

\medskip 
\noindent
We remark that the point $(x:y:z)=(1:1:1)$ always corresponds with $(X:Y:Z) =$ $(1:1:1)$.
For $X=Y=Z=1$ we have the symmetric matrix:
\begin{align}\label{eq:111}
& M(\mu)=\left( \begin{array}{ccc} 1 & 3\mu & 3\mu \\
                          3\mu & 1 & 3\mu \\
                           3\mu & 3\mu & 1
         \end{array} \right), \ \ \ \mu \in  (-1/6,1/2]
\end{align}

\noindent
which is positive definite for $\mu \in (-1/6,1/3)$ and indefinite for 
$\mu \in (1/3,1/2]$. 

\medskip 
We should remain aware of the fact that $\rho=0$ and $\rho=-1/27$ in (\ref{eq:Det})
give singular cubics, hence the framework ${\cal F}(\lambda)$ will be {\em non-regular} for $\lambda = 0, 1/3, 1/2, 2/3, 1$.

\medskip\noindent{\bf Particular cases.} We consider the three cases $\lambda=1/6,1/3,5/12$ and follow the algorithm.

\medskip
For $\lambda=1/6$, by STEP 3, we obtain the invariants: 

$$ S=-\frac{2287}{4000752},\ T=-\frac{2021723}{18525482136},\ 
\Delta=\frac{1000000}{22067482534159923} $$

\medskip \noindent
In STEP 4 we find $k=25.6407$ and, via \eqref{eq:111}, the output is {\em YES: the regular framework admits an auxetic deformation}.

\medskip 
For $\lambda=1/3$, by STEP 3, we obtain the invariants:

$$ S=-\frac{1}{11664}, \ T=\frac{1}{157464},\ \Delta=0 $$

\medskip \noindent
Since $\Delta=0$ the algorithm stops with the message that {\em the framework is NOT REGULAR.}

\medskip 
For $\lambda=5/12$, by STEP 3, we obtain the invariants:

$$ S=-\frac{9973}{25625808}, \ T=-\frac{45441143}{760048652376}, \ 
\Delta=\frac{7353062500}{37144672966729275363}   $$

\medskip \noindent
In STEP 4 we find $k=10.6042$ and, via \eqref{eq:111}, the final output is {\em NO: the regular framework does not admit an auxetic deformation}.

\medskip
\noindent
{\bf Conclusion.}
Using ideas from the theory of elliptic curves and invariants of ternary cubics, we have proposed an algorithm  for deciding if a regular periodic bar-and-joint framework in 3D, with three degrees of freedom, admits an auxetic deformation. 

\medskip
\noindent
{\bf Acknowledgements.} The authors are grateful for the hospitality of the Institute for Computational and Experimental Research in Mathematics (ICERM) at Brown University during the Fall 2016, when this work was completed.


\begin{thebibliography}{10}

\bibitem{bonifant:milnor:cubicCurves:arxiv:2016}
A.~Bonifant and John Milnor.
\newblock Smooth cubic curves over {$C$} and {$R$}.
\newblock 2016.
\newblock arXiv: 1603.09018.

\bibitem{borcea:streinu:periodicFlexibility:2010}
Ciprian~S. Borcea and Ileana Streinu.
\newblock Periodic frameworks and flexibility.
\newblock {\em Proceedings of the Royal Society A}, 466(2121):2633--2649, 2010.

\bibitem{borcea:streinu:minimallyRigidPeriodic:BLMS:2011}
Ciprian~S. Borcea and Ileana Streinu.
\newblock Minimally rigid periodic graphs.
\newblock {\em Bulletin of the London Mathematical Society}, 43:1093--1103,
  2011.

\bibitem{borcea:streinu:frameworksCrystallographicSymmetry:philTransAndArxiv:2013}
Ciprian~S. Borcea and Ileana Streinu.
\newblock Frameworks with crystallographic symmetry.
\newblock {\em Philosophical Transactions of the Royal Society of London Ser.
  A:}, 372(20120143), 2013.

\bibitem{borcea:streinu:kinematicsExpansive:ark14:2014}
Ciprian~S. Borcea and Ileana Streinu.
\newblock Kinematics of expansive planar periodic mechanisms.
\newblock In {\em Adv. Robot Kinematics (ARK'14)}, pages 395--408. Springer,
  2014.

\bibitem{borcea:streinu:geometricAuxetics:RSPA:arxiv:2015}
Ciprian~S. Borcea and Ileana Streinu.
\newblock Geometric auxetics.
\newblock {\em Proceedings of the Royal Society A}, 471(20150033), 2015.
\newblock arXiv:1501.03550.

\bibitem{borcea:streinu:liftingsStresses:dcg:arxiv:2015}
Ciprian~S. Borcea and Ileana Streinu.
\newblock Liftings and stresses for planar periodic frameworks.
\newblock {\em Discrete and Computational Geometry}, 53(4):747--782, 2015.
\newblock arxiv:1501.03549.

\bibitem{borcea:streinu:NewPrinciplesAuxeticDesign:arxiv:2016}
Ciprian~S. Borcea and Ileana Streinu.
\newblock New principles for auxetic periodic design.
\newblock Submitted, 2016.
\newblock arxiv:1608.02104.

\bibitem{chen:wang:inflectionPoints:2003}
Falai Chen and Wenping Wang.
\newblock Computing real inflection points of cubic algebraic curves.
\newblock {\em Computer Aided Geometric Design}, 20:101--117, 2003.

\bibitem{dolgachev:classicalAlgGeom:2012}
Igor~V. Dolgachev.
\newblock {\em Classical algebraic Geometry}.
\newblock Cambridge University Press, 2012.

\bibitem{dove:displacive:1997}
Martin~T. Dove.
\newblock Theory of displacive phase transitions in minerals.
\newblock {\em American Mineralogist}, 82:213--244, 1997.

\bibitem{elipe:lantada:auxeticGeometries:2012}
J.~C.~A. Elipe and A.~D. Lantada.
\newblock Comparative study of auxetic geometries by means of computer-aided
  design and engineering.
\newblock {\em Smart Materials and Structures}, 21:105004, 2012.

\bibitem{evans:etAl:molecularNetwork:Nature:1991}
Kenneth~E. Evans, M.~A. Nkansah, I.~J. Hutchinson, and S.~C. Rogers.
\newblock Molecular network design.
\newblock {\em Nature}, 353:124--125, 1991.

\bibitem{greaves:surveyPoissonRatios:resNotesRoyalSoc:2013}
G.~N. Greaves.
\newblock Poisson's ratio over two centuries: challenging hypotheses.
\newblock {\em Notes and Records of the Royal Society of London}, 67:37--58,
  2013.

\bibitem{greaves:lakes:etAl:PoissonRatio:2011}
G.~N. Greaves, A.~I. Greer, R.~Lakes, and T.~Rouxel.
\newblock Poisson's ratio and modern materials.
\newblock {\em Nature Materials}, 10:823--837, 2011.

\bibitem{grima:etAl:doZeolitesNegativePoission:advMaterials:2000}
Joseph Grima, Rosie Jackson, Andrew Alderson, and Kenneth~E. Evans.
\newblock Do zeolites have negative {P}oisson's ratios?
\newblock {\em Advanced Materials}, 12(24):1912--1918, 2000.

\bibitem{huang:chen:negativePoisson:2016}
C.~Huang and L.~Chen.
\newblock Negative {P}oisson's ratio in modern functional materials.
\newblock {\em Advanced Materials}, 28:8079--8096, 2016.

\bibitem{lakes:negativePoisson:1987}
R.~Lakes.
\newblock Foam structures with a negative {P}oisson's ratio.
\newblock {\em Science}, 235:1038--1040, 1987.

\bibitem{lee:singer:thomas:microNanoMaterials:advMat:2012}
Jae-Hwang Lee, Jonathan~P. Singer, and Edwin~L. Thomas.
\newblock Micro-/nanostructured mechanical metamaterials.
\newblock {\em Advanced Materials}, 24:4782--4810, 2012.

\bibitem{lemon:manchoso:ye:LowRankSDP:foundTrends:2016}
Alex Lemon, Anthony Man-Cho So, and Yinyu Ye.
\newblock Low-rank semidefinite programming: Theory and applications.
\newblock {\em Foundations and Trends in Optimization}, 2(1-2):1--156, 2016.

\bibitem{megaw:crystalStructures:1973}
H.~D. Megaw.
\newblock {\em Crystal structures: a working approach}.
\newblock W.B.Saunders Co., 1973.

\bibitem{mitschkeEtAl:auxetic:rspa:2013}
H.~Mitschke, V.~Robins, K.~Mecke, and G.~E. Schr\"{o}der-Turk.
\newblock Finite auxetic deformations of plane tessellations.
\newblock {\em Proceedings of the Royal Society A}, 469(20120465), 2013.

\bibitem{ottem:etAl:quarticSpectrahedra:2015}
John~Christian Ottem, Kristian Ranestad, Bernd Sturmfels, and Cynthia Vinzant.
\newblock Quartic spectrahedra.
\newblock {\em Mathematical Programming Ser.B}, 151(2):585--612, 2015.
\newblock arXiv:1311.3675.

\bibitem{porkolab:khachian:semidefinite:1997}
L.~Porkolab and Leonid Khachian.
\newblock On the complexity of semidefinite programs.
\newblock {\em Journal of Global Optimization}, 10(4):351--365, 1997.

\bibitem{sottile:ternaryQuartics:2004}
Victoria Powers, Bruce Reznick, Claus Scheiderer, and Frank Sottile.
\newblock A new approach to {H}ilbert's theorem on ternary quartics.
\newblock {\em C. R. Acad. Sci. Paris, Ser. I}, 339:617--620, 2004.

\bibitem{reisEtAl:designerMatter:2015}
P.~M. Reis, H.~M. Jaeger, and M.~van Hecke.
\newblock Designer matter: A perspective.
\newblock {\em Extreme Mechanics Letters}, 5:25--29, 2015.

\bibitem{siddornEtAl:systematicTopology:2015}
M.~Siddorn, F.~X. Coudert, Kenneth~E. Evans, and Arnaud Marmier.
\newblock A systematic typology for negative {P}oisson's ratio materials and
  the prediction of complete auxeticity in pure silica zeolite {JST}.
\newblock {\em Physical Chemistry Chemical Physics}, 17:17927--17933, 2015.

\bibitem{sturmfels:algorithmsInvariantTheory:2008}
Bernd Sturmfels.
\newblock {\em Algorithms in Invariant Theory}.
\newblock Springer, 2008.

\bibitem{vinzant:whatIsSpectrahedron:2014}
Cynthia Vinzant.
\newblock What is ... a spectrahedron?
\newblock {\em Notices of the American Mathematical Society}, 61:492--494,
  2014.

\bibitem{yeganeh:etAl:elasticityCristobalite:1992}
Amir Yeganeh-Haeri, Donald~J. Weidner, and John Parise.
\newblock Elasticity of alpha-cristobalite: A silicon dioxide with a negative
  {P}oisson's ratio.
\newblock {\em Science}, 257(5070):650--652, July 1992.

\end{thebibliography}
\end{document}